\documentclass[12pt]{amsart}
\usepackage{amsmath,amsthm,amsfonts,amssymb,color,tikz}
\begin{document} 
\newcommand{\XXX}{\message{XXX}.{\bf XXX !}}
\newcommand{\changed}[1]{\textcolor{red}{#1}}
\newcommand{\anfangchanged}{\color{red}}
\newcommand{\finchanged}{\color{black}}
\newcommand{\modulo}[1]{\text{ mod }#1}
\newcommand{\A}{{\mathbb A}}
\newcommand{\B}{{\mathbb B}}
\newcommand{\C}{{\mathbb C}}
\renewcommand{\H}{{\mathrm H^+}}

\newcommand{\F}{{\mathbb F}}
\newcommand{\N}{{\mathbb N}}
\newcommand{\Q}{{\mathbb Q}}
\newcommand{\Z}{{\mathbb Z}}
\renewcommand{\P}{{\mathbb P}}
\newcommand{\R}{{\mathbb R}}
\newcommand{\rc}{\subset}
\newcommand{\rank}{\mathop{rank}}
\newcommand{\trace}{\mathop{tr}}
\newcommand{\dimc}{\mathop{dim}_{\C}}
\newcommand{\Lie}{\mathop{Lie}}
\newcommand{\Spec}{\mathop{Spec}}
\newcommand{\Auto}{\mathop{{\rm Aut}_{\mathcal O}}}
\newcommand{\alg}[1]{{\mathbf #1}}
\newcommand{\tensor}{\otimes}
\newtheorem{lemma}{Lemma}[section]
\newtheorem*{definition}{Definition}
\newtheorem*{claim}{Claim}
\newtheorem{corollary}{Corollary}
\newtheorem*{Conjecture}{Conjecture}
\newtheorem*{SpecAss}{Special Assumptions}
\newtheorem{example}{Example}
\newtheorem*{remark}{Remark}
\newtheorem*{observation}{Observation}
\newtheorem*{fact}{Fact}
\newtheorem*{remarks}{Remarks}
\newtheorem{proposition}[lemma]{Proposition}
\newtheorem{theorem}[lemma]{Theorem}
\numberwithin{equation}{section}
\def\labelenumi{\rm(\roman{enumi})}
\title{%
Unexpected Examples for the Kobayashi Pseudo distance
}
\author {J\"org Winkelmann}
\begin{abstract}
  We present some unexpected examples related to the
  Kobayashi pseudodistance. There are unramified coverings where
  the Kobayashi pseudo distance vanishes on the base manifold
  but not on the total space. The Kobayashi pseudodistance may
  vanish identically even if the Kobayashi-Royden Pseudo
  metric does not.
\end{abstract}
\keywords{Kobayashi pseudodistance, covering}
\subjclass{Complex Variables
32Q45}%
%
\address{%
J\"org Winkelmann \\
IB 3-111\\
Lehrstuhl Analysis II \\
Fakult\"at f\"ur Mathematik \\
Ruhr-Universit\"at Bochum\\
44780 Bochum \\
Germany
}
\email{joerg.winkelmann@rub.de\newline
    ORCID: 0000-0002-1781-5842
}
\thanks{
  {\em Acknowledgement.}
  The research for this article was initiated by a question
  raised by J.~Aryanpekhar at a workshop organized by D.~Brotbek in Nancy.
}
\maketitle

\section{Introduction}

\subsection{Kobayashi Pseudo Distance}

For every complex space there exists an intrinsically defined
pseudo-distance, the {``\em Kobayashi pseudo distance''}.
See \cite{Kobayashi}. It plays an important r\^ole in complex
geometry.

From the definition of the Kobayashi pseudodistance via disc chains
one may derive the following {``\em covering formula''}
(see \cite{Kobayashi}, Theorem. 3.2.8)

\begin{proposition}
  Let $\pi:X\to Y$ be an unramified covering of complex manifolds, $p,q\in Y$,
  $\tilde p\in\pi^{-1}(p)$.
  Then
  \begin{equation}\label{cover-formel}
    d_Y(p,q)=\inf_{\tilde q\in\pi^{-1}(q)}d_X(\tilde p,\tilde q)
  \end{equation}
\end{proposition}

With the aid of this formula one may deduce:
{\em Given an unramified covering $X\to Y$, the total
  space $X$ is hyperbolic (in the sense of Kobayashi)
  \footnote{A complex space is {\em (Kobayashi) hyperbolic} iff its
    Kobayashi pseudo distance is a distance.}
  if and only
  if the base space  $Y$ is hyperbolic}.

This raises the question whether the same is also true for total degeneration
of the Kobayashi pseudo distance, i.e., given an unramified
covering $\pi:X\to Y$, is the condition $d_X\equiv 0$ equivalent to
$d_Y\equiv 0$?

The above stated covering formula \eqref{cover-formel} immediately
yields
the implication $d_X\equiv 0\implies d_Y\equiv 0$

But for the opposite direction
we construct counter examples: There are infinite unramified
coverings $\pi:X\to Y$ with $d_Y\equiv 0$, but $d_X\not\equiv 0$.
(Theorem~\ref{main}, Theorem~\ref{main2}).

However, for finite coverings, there 
is a positive answer (Proposition~\ref{finite}).

\subsection{Kobayashis question}

In \cite{K70} (p.~48), S.~Kobayashi raised the question whether the infimum
in \eqref{cover-formel} is always attained. In our examples
as described in Theorem~\ref{main} and Theorem~\ref{main2} we have
  \begin{equation}
    0=\inf_{\tilde q\in\pi^{-1}(q)}d_X(\tilde p,\tilde q)
  \end{equation}
  with $d_X(\tilde p,\tilde q)>0$ for all
  $\tilde q\in\pi^{-1}(q)$. Thus the infimum is not attained.

  However, such examples have been known before.
  Namely, in \cite{Z} Zwonek constructed Reinhardt domains with such
  properties.
  Here the covering map goes from
  $\Delta^ *\times\C$ (with $\Delta^ *=\{z\in\C:0<|z|<1\}$)
  to the Reinhard domain
  \[
  D_\alpha=\{(z_1,z_2)\in\C^ *\times\C^*:|z_1|^ {\alpha_1}|z_2|^ {\alpha_2}<1\}
  \quad\quad (\alpha_1<0<\alpha_2,\ \ \frac{\alpha_1}{\alpha_2}\not\in\Q)
  \]
  
  Other examples are provided by the universal
  covering of Inoue surfaces (\cite{In}).
  These are non-K\"ahler compact surfaces whose universal covering
  is $\Delta\times\C$.
  
  Furthermore in \cite{WK} we discussed a non-compact surface
  which arises as quotient of $\Delta\times\C$ by an action
  of $\Z^ 2$.

  In all these examples we have a covering where both on the base
  and on the total space the Kobayashi pseudo distance neither
  vanishes nor is a proper distance.

  The total space is in all three examples
  $\Delta\times\C$ or $\Delta^ *\times\C$ and the covering map
  maps the complex lines in $\Delta\times\C$
  resp.~$\Delta^*\times\C$
  to the base manifold
  in such a way that the closure is real three-dimensional.
  Then, if $p,q$ are not in the image of one complex line, but
  both in the closure of the same complex line, the infimum
  in \eqref{cover-formel} is not attained.
  
  But these examples do not answer
  whether
  the implication
  \[
  d_X\equiv 0\implies d_Y\equiv 0
  \]
  holds, because the Kobayashi pseudodistance on the base
  manifold does not vanish identically.
  
\subsection{The Infinitesimal Kobayashi-Royden pseudometric}

Let $X$ be a complex manifold, $x\in X$, $v\in T_xX$.

Then the {\em Kobayashi Royden pseudometric} $F_X$
(see \cite{Royden}) is
  defined as
  \[
  F_X(v)=\inf\{ \lambda\in\R^+:\exists f:\Delta\to X \text{ holo. }
  f(0)=x, \lambda f'(0)= v
  \}
  \]

  Royden proved that $F_X$ is semi-continuous and therefore
  locally $L^1$ and demonstrated that $d_X(x,y)$ can be realized as
  \[
  d_X(x,y)=\inf_\gamma \int^1_0 F_X(\gamma'(t))dt
  \]
  with the infimum taken over all integration paths $\gamma:[0,1]\to X$ with $\gamma(0)=x$
  and $\gamma(1)=y$.
 (See \cite{Royden},\cite{Kobayashi}).
  
  From this it follows that $d_X$ vanishes if $F_X$ vanishes.

  We give an example of a complex manifold with vanishing $d_X$
  but non-vanishing $F_X$ (Theorem~\ref{main-inf}).

  \subsection{Fiber bundles}
  We show that there is a locally trivial holomorphic fiber bundle
  such that the fiber is Kobayashi hyperbolic, but the total space has
  vanishing Kobayashi pseudo distance (Theorem~\ref{bundle}).

  It is an interesting point that this bundle is {\em holomorphically}
  locally trivial, since there are easier examples which are only
  differentiably trivial. For example, let
  \[
  E=\{[x:y:z]\in\P_2(\C): |x|^2+|y|^ 2>|z|^2\}
  \]
  Then $E\ni[x,y,z]\mapsto[x:y]\in\P_1(\C)$ defines
  a $C^ \infty$-locally trivial bundle with fiber isomorphic to $\Delta$
  and $d_E\equiv 0$. But this bundle is not
  holomorphically locally trivial.

  \section{Preparation}

\subsection{Lie Subgroups of $SL_2(\R)$}

\begin{lemma}\label{sub-solve}
  Let $H$ be a closed subgroup of $G=SL_2(\R)$.

  If $0<\dim H<3$, then $H$ is solvable.
\end{lemma}

\begin{proof}
  Let $H^0$ denote the connected component of $H$ containing the
  neutral element. This is a connected Lie subgroup. Since connected
  Lie subgroups of a given Lie group correspond to the Lie subalgebras
  of the corresponding Lie algebra, $H^0$ is conjugate to one
  of the following:
  \[
  B=
  \left\{ \begin{pmatrix} \lambda & c \\ 0 & \lambda^{-1} \\
  \end{pmatrix}
  : \lambda\in\R^*, c\in \R
  \right\},
  \quad
  U=
  \left\{ \begin{pmatrix} 1 & c \\ 0 & 1 \\
  \end{pmatrix}
  : c\in \R
  \right\},
  \]
   
  \[
  K=
  \left\{ \begin{pmatrix} \cos\theta & \sin\theta \\
    -\sin\theta & \cos\theta \\
  \end{pmatrix} \right\}
  \text{ or }
  T^+=
  \left\{ \begin{pmatrix} \lambda &  0\\ 0 & \lambda^{-1} \\
  \end{pmatrix}
  : \lambda\in\R^+
  \right\}
  \]
  Observe that conjugation by any element of $H$ stabilizes $H^ 0$.
  Therefore $H$ is contained in
  \[
  N_G(H^ 0)=\{g\in G: gH^ 0g^ {-1}=H^ 0\}
  \]
  A case-by-case check verifies that $N_G(H^0)$ is always solvable,
  if $H^0$ is a non-trivial (i.e. $H\ne G$) connected Lie subgroup
  of $SL_2(\R)$:
  $H^ 0$ is conjugate to one of the above listed subgroups $B$, $U$ and
  $T^ +$, and
  \[
  N_G(B)=B,\quad N_G(U)=B,
  \]
  \[
  N_G(T^+)=
  \left\{ \begin{pmatrix} \lambda &  0\\ 0 & \lambda^{-1} \\
  \end{pmatrix}
  : \lambda\in\R^*
  \right\}
  \cup
  \left\{
  \begin{pmatrix} 0 & \lambda \\ -\lambda^ {-1} & 0 \\
  \end{pmatrix}:\lambda\in\R^*
  \right \}
  \]
  and all these groups are easily seen to be solvable.
\end{proof}

\begin{lemma}\label{common-eigen}
  Let $H$ be a closed subgroup of $SL_2(\R)$.
  Let $g,h\in H$.

  Assume $0<\dim(H)<3$.
  Then $g^2$ and $h^2$ have a common eigen vector.
\end{lemma}

\begin{proof}
  By the above considerations $H$ has at most two connected component.
  Thus $g^2,h^2\in H^0$ where $H^0$ denotes the connected
  component of $H$ which contains the neutral element.
  $H^0$ is connected and solvable (Lemma~\ref{sub-solve}).
  Hence there is a common eigen vector.
\end{proof}

\subsection{Dense subgroups}

The existence of dense subgroups with $2$ generators
in $SL_2(\R)$ follows from
general results of \cite{BG}, but may also be shown explicitly.

\begin{lemma}
  Let
  \[
  A=\begin{pmatrix} 2 & 1 \\ 0 & \frac 12
  \end{pmatrix},\quad
  B=\begin{pmatrix} 3 & 0 \\ 1 & \frac 13\\
  \end{pmatrix}
  \]
  Then $A$ and $B$ generate a dense subgroup of $SL(2,\R)$.
\end{lemma}

\begin{proof}

  We calculate:
  \[
  A^n=\begin{pmatrix}
  2^n & \frac{4^n-1}{2^{n-1}3}\\0 & 2^{-n}\\
  \end{pmatrix}
  ,\quad
  B^m=
  \begin{pmatrix}
    3^m & 0 & \\
    \frac{9^m-1}{3^{m-1}8} & 3^{-m}\\
  \end{pmatrix}
  \]
  
Eigen vectors:
\begin{align*}
A\begin{pmatrix} 1 \\ 0\\
\end{pmatrix}
&= 2\begin{pmatrix} 1 \\ 0\\
\end{pmatrix}\\
A\begin{pmatrix} 2 \\ -3\\
\end{pmatrix}
&=  \frac 12\begin{pmatrix} 2 \\ -3\\
\end{pmatrix}\\
B\begin{pmatrix} 0 \\ 1\\
\end{pmatrix}
&= \frac 13\begin{pmatrix} 0 \\1 \\
\end{pmatrix}\\
B\begin{pmatrix} 8 \\ 3\\
\end{pmatrix}
&= 3\begin{pmatrix} 8 \\ 3\\
\end{pmatrix}\\
\end{align*}
Let $H$ denote the closure of the group generated by $A$ and $B$.
Due to Cartans theorem this is a Lie subgroup \cite{Cartan}.

Assume first $0<\dim H<3$. The above calculations show that
there are no $n,m\in\N$ such that $A^n$ and $B^m$ have a common eigen vector.
This contradicts lemma~\ref{common-eigen}. Hence this case is impossible.

Now assume that $H$ is discrete.

Observe that $(m,n)\mapsto 2^n3^m$ is an injective map from
$\Z\times\Z$ to $\R^+$ and that the image is dense.
In particular, $1$ is an accumulation point of
\[
\{2^n3^m:n,m\in\Z\}
\]
We consider $A^nB^m$.
\[
A^nB^m=\begin{pmatrix}
3^n2^m & 3^n\frac{2^{2m}-1}{2^{n-1}3} \\
2^m \frac{3^{2n}-1}{3^n8}
& \frac{(3^{2n-1})(2^{2m}-1)}{3^n2^{m-1}24}+3^{-n}2^{-m}\\
\end{pmatrix}
\]
We choose a sequence among the $(m,n)\in\Z^2$ such that
$m\to\infty$, $n\to\infty$ and $2^n3^m\to 1$.
Then
\[
A^nB^m \to
\underbrace{\begin{pmatrix}
  1 & \frac 23 \\
  \frac 18 & \frac{13}{12} \\
\end{pmatrix}}_C
\]
Since $2^n3^m$ equals the upper left entry of the $2\times 2$-matrix
$A^nB^m$ and $2^n3^m\to 1$, but $2^n3^m\ne 1\forall (n,m)\ne(0,0)$,
we conclude that $C$ is an accumulation point of
the set $\{A^nB^m:n,m\in\Z\}$. Thus
the subgroup $H$ ( which contains $A$ and $B$) can not be discrete.

It follows that $H=SL_2(\R)$, since the case $0<\dim(H)<3$ was
already ruled out.
\end{proof}

\subsection{Free discrete subgroups of $PSL_2(\R)$}

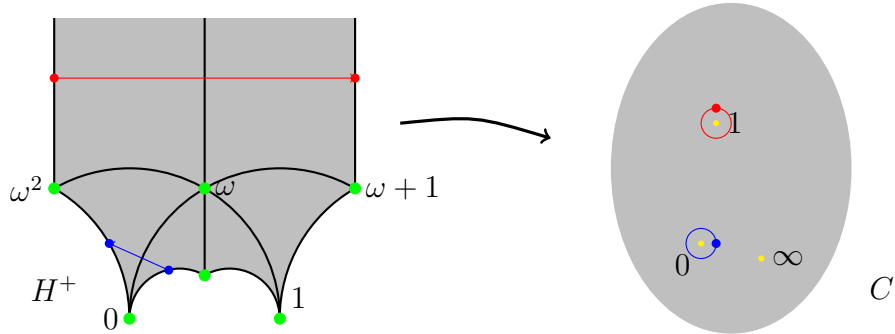
\begin{center}
  \begin{tikzpicture}[scale=2]
    \fill [lightgray] (-0.5,2) -- (-0.5,0.866)
    arc [start angle=60,end angle=0,radius=1]
     arc [start angle=180, end angle=60, radius=0.333]
    arc [start angle=120, end angle=0, radius=0.333]
    arc [start angle=180,end angle=120,radius=1]
    -- (1.5,0.866)  -- (1.5,2) -- (-0.5,2);
    \draw [black,thick] (-0.5,2) -- (-0.5,0.866)
    arc [start angle=60,end angle=0,radius=1]
     arc [start angle=180, end angle=60, radius=0.333]
    arc [start angle=120, end angle=0, radius=0.333]
    arc [start angle=180,end angle=120,radius=1]
    -- (1.5,0.866)  -- (1.5,2);
    \draw [black,thick] (1,0) arc [start angle=0,end angle=120,radius=1];
    \draw [black,thick] (0,0) arc [start angle=180,end angle=60,radius=1];
    \draw [black,thick]  (0.5,0.288) --  (0.5,2);
    \draw [red,->] (-0.5,1.6) -- (1.5,1.6);
    \draw [blue,->] (0.2617,0.32256) -- (-0.134,0.5);
    \node at (-0.5,0.2) {$H^+$};
    \fill [green] (-0.5,0.866) node [left,black] {$\omega^2$} circle ( 0.04);
    \fill [green] (1.5,0.866) node [right,black] {$\omega+1$} circle ( 0.04);
    \fill [green] (0.5,0.866) node [right,black] {$\omega$} circle ( 0.04);
    \fill [green] (0.5,0.288) circle ( 0.04);
    \fill [green] (0,0) node [black,left] {$0$} circle ( 0.04);
    \fill [green] (1,0) node [black,above right] {$1$} circle ( 0.04);
    \fill [blue] (-0.134,0.5) circle ( 0.03);
    \fill [blue]  (0.2617,0.32256)  circle ( 0.03);
    \fill [red] (-0.5,1.6)  circle ( 0.03);
    \fill [red] (1.5,1.6)   circle ( 0.03);
    \fill [lightgray] (4,1) ellipse [x radius=0.8,y radius=1.1];
    \fill [yellow] (3.8,0.5) circle (0.02);
    \fill [yellow] (4.2,0.4) circle (0.02) node [black,right] {$\infty$};
    \fill [yellow] (3.9,1.3) circle (0.02);
    \draw [blue]  (3.8,0.5) circle (0.1) node [black,below left] {$0$};
    \draw [red]  (3.9,1.3)  circle (0.1)  node [black,right] {$1$};
    \fill [blue] (3.9,0.5) circle (0.03);
    \fill [red] (3.9,1.4) circle (0.03);
    \node at (5.0,0.2) {$C$};
    \draw [very thick,black,->] (1.8,1.3) .. controls (2.3,1.35) ..
    (2.8,1.2);
  \end{tikzpicture}
\end{center}
Due to \cite{JW1} every semisimple Lie group
contains both discrete and non-discrete
free subgroups with $2$ generators;
for the case $SL_2(\R)$, see also \cite{JW2}.

However, from the theory of {\em modular curves}
there are also well-known explicit examples.

Consider
\[
\Gamma(2)=\{ A\in PSL_2(\Z): A\equiv I \text{ \sl mod }\ 2\}.
\]

This is a discrete and torsion-free subgroup of $PSL_2(\R)$
which acts properly discontinuously on the upper half plane
$\H=\{z\in\C:\Im(z)>0\}$ with the fundamental domain shown
in the diagram.

The corresponding {\em  modular curve} $C=\H/\Gamma(2)$
may be compactified by adding one point for every element
of
\[
\P_1(\Q)/\Gamma(2)\simeq \P_1(\F_2) \simeq \{0,1,\infty\}.
\]
Then $C$ equals $\P_1$ with three points removed.
Since $\Gamma(2)$ is torsion-free,
we obtain that
\[
\Gamma(2)\simeq \pi_1(C)\simeq\pi_1\left(
\P_1\setminus\{0,1,\infty\}\right)
\]
is a free group with two generators.

More explicitly, let
\[
U=\begin{pmatrix} 1 & 2 \\
0 & 1 \\
\end{pmatrix},\quad
V=\begin{pmatrix} 1 & 0 \\
2 & 1 \\
\end{pmatrix}
\]

Then $U$ yields a circle around $\infty$ in $C$, denoted red in the
diagram and $V$ yields a circle around $0$ in $C$,
painted blue.

Thus $\Gamma(2)$ is a free discrete subgroup of $PSL_2(\R)$
generated by the two elements $U$ and $V$.

Since $\Gamma(2)$ is free, the natural embedding
of $\Gamma(2)$ into $PSL_2(\R)$ lifts to an
embedding into $SL_2(\R)$
by choosing an arbitrary lift for its two generators.

\subsection{The key lemma}

We summarize now the results of the above considerations
in the form which we will be used subsequently.

\begin{proposition}\label{emb-f2}
  Let $F_2$ be the free group with two generators. Then:

  \begin{itemize}
  \item
    There is an embedding $j_1$  of    $F_2$  into $SL_2(\R)$
    as a discrete subgroup.
  \item
    There is a group homomorphism $j_2$  of    $F_2$
    onto a dense subgroup of $PSL_2(\R)$.
    \footnote{In fact $j_2$ may be chosen to be injective.}
  \end{itemize}

  Explicitly we may choose $j_1$ by mapping the two generators of $F_2$
  to the elements of $PSL_2(\R)$ represented by the matrices
  \[
  \begin{pmatrix} 2 & 1 \\ 0 & \frac 12
  \end{pmatrix},\quad
  \begin{pmatrix} 3 & 0 \\ 1 & \frac 13\\
  \end{pmatrix}
  \]
  and $j_2$ by mapping the generators of $F_2$ to
  \[
  \begin{pmatrix} 1 & 2 \\
  0 & 1 \\
  \end{pmatrix},\quad
  \begin{pmatrix} 1 & 0 \\
  2 & 1 \\
  \end{pmatrix}
  \]
  respectively.
\end{proposition}

\section{Main construction}

\begin{theorem}\label{main}
  Let $X=\H\times SL_2(\C)$ with $\H=\{z\in\C:\Im(z)>0\}$.
  There exists a properly discontinuous and free action $\mu$ of
  the free group $F_2$ with $2$ generators on $X$ such that
  $d_Y\equiv 0$ for $Y=X/\mu$.
\end{theorem}

\begin{remark}
  Note that $d_X\not\equiv 0$, since $\H$ is hyperbolic. In fact,
  it is well-known that
  \[
  d_{\H}(z,w)=arcosh\left( 1 + \frac{|z-w|^2}{2\Im(z)\Im(w)}\right)
  \]
  Hence this gives an example of an unramified covering $X\to Y$ with
  $d_Y\equiv 0$, but $d_X\not\equiv 0$.
\end{remark}

\begin{proof}
  Using Proposition~\ref{emb-f2} we choose embeddings $j_1:F_2\to SL_2(\R)$
  and $j_2:F_2\to SL_2(\R)\subset SL_2(\C)$ such that $j_1$ has dense
  image in $SL_2(\R)$ and $j_2(F_2)$ is discrete.
  The action $\mu$ of $F_2$ on $X$ is now defined by
  \[
  \mu(\alpha):(z,A)\mapsto\left( \phi(j_1(\alpha))(z),j_2(\alpha)\cdot A\right)
  \]
  where $\phi$ is the usual action by M\"obius transformations, i.e.,
  \[
  \phi\begin{pmatrix} a & b \\ c & d \\
  \end{pmatrix}(z)
  =
  \frac{az+b}{cz+d},
  \]
  and $j_2(\alpha)\cdot A$ denotes the product with respect to the
  group law on $SL_2(\C)$.
  Discreteness of $j_2(F_2)$ implies that this action $\mu$ is free and
  properly discontinuous. Thus we obtain an unramified covering
  $\pi:X\to X/\mu\stackrel{def}{=}Y$.

  Let $(x,A),(y,B)\in \H\times SL_2(\C)=X$ be arbitrary.
  Since $j_1(F_2)$ is dense and since $SL_2(\R)$ acts transitively
  on $\H$, there is a sequence $\alpha_n\in F_2$ with
  \[
  \lim_{n\to\infty}\left( \phi(j_1(\alpha_n))(y)\right)=x
  \]
  Now
  \[
  d_X\left( (x,A),(y,B)\right)=d_{\H}(x,y),
  \]
  because $SL_2(\C)$ is a complex Lie group and therefore has
  vanishing Kobayashi pseudodistance
  (\cite{Kobayashi}, Ex. 3.1.22).
  By the covering formula \eqref{cover-formel} we have
  \begin{align*}
  d_Y(\pi(x,A),\pi(y,B))
  &=\inf_{\alpha\in F_2}
  d_X((x,A),\mu(\alpha)(y,B))\\
&  =\inf_{\alpha\in F_2}
  d_{\H}(x,\phi(j_1(\alpha))(y))\\
  \end{align*}
  But
  \[
 \inf_{\alpha\in F_2}
 d_{\H}(x,\phi(j_1(\alpha))(y))=0
 \]
 because
  \[
  \lim_{n\to\infty}\left( \phi(j_1(\alpha_n))(y)\right)=x.
  \]
  Therefore
  \[
  d_Y\left(\pi(x,A),\pi(y,B)\right)=0
  \]
  for {\em every} $(x,A)$ and $(y,B)$ in $X$.
  Thus $d_Y\equiv 0$.
\end{proof}

\section{Another example, based on the construction of Margulis}

Semisimple Lie groups like $SL_2(\C)$ always carry non-trivial topology.
Thus one might ask, whether there is an example of
an unramified cover as in Theorem~\ref{main} above
with a contractible total space.

This is indeed the case. We recall that Margulis constructed
a properly discontinuous free action  of the free group
with two generators $F_2$ on $\R^3$
by affin-linear transformations
(\cite{Margulis}).

If in the preceding construction we replace $SL_2(\C)$ by $\C^3$
and the action of $F_2$ via $j_2$ and right multiplication by the
complexified action of Margulis, we obtain an unramified covering
$X'\to Y'$ with $X'=\H\times\C^3$ and $d_{Y'}\equiv 0$.
In this way we establish the below result.

\begin{theorem}\label{main2}
  There exists a complex manifold $Y'$ with identically vanishing
  Kobayashi pseudo distance $d_{Y'}$ whose universal covering is
  biholomorphic to
  \[
  \H\times\C^ 3=\{z\in\C^ 4: \Im(z_1)>0\}.
  \]
\end{theorem}

 \section{Infinitesimal Pseudodistance}

 Let $\pi:X\to Y$ be an unramified covering, $x\in X$, $y=\pi(x)$.
 Every holomorphic
 map $f$ from the unit disc $\Delta$ to $X$ mapping $0$ to $x$
 yields a map $g$ from the unit disc to $Y$ by concatenation
 with $\pi$ (i.e., $g=\pi\circ f$). Conversely, since $\Delta$ is
 simply-connected, every holomorphic map from $Y$ lifts to a map to $X$.
 
 Therefore, the definition of the Kobayashi Royden pseudodistance implies
 that $\pi^*F_Y=F_X$, i.e.,
 \[
 F_Y\left((D\pi)(v)\right)=F_X(v) \quad \forall v\in T_xX
 \]

 Thus, given an unramified covering $\pi:X\to Y$, we have $F_X\equiv 0$
 if and only if $F_Y\equiv 0$.

 In the preceding section we constructed unramified coverings
 $\pi:X\to Y$ where the Kobayashi pseudo distance of $Y$ is vanishing
 identically, but $d_X\not\equiv 0$.

 Since
  \[
  d_X(x,y)=\inf_\gamma \int^1_0 F_X(\gamma'(t))dt
  \]
  where the infimum is taken over all paths
  $\gamma:[0,1]\to X$ with $\gamma(0)=x$
  and $\gamma(1)=y$
  (see \cite{Royden}), the non-vanishing of $d_X$ implies
  that $F_X$ is not vanishing.
  The non-vanishing of $F_X$ in turn implies that $F_Y$ is not
  vanishing.

  Thus every unramified covering $\pi:X\to Y$ with $d_Y\equiv 0$
  and $d_X\not\equiv 0$ yields an example of a complex manifold $Y$
  with $d_Y\equiv 0$ but $F_Y\not\equiv 0$.
  
  In combination with Theorem~\ref{main} or Theorem~\ref{main2}
  this yields:

  \begin{theorem}\label{main-inf}
    There exists a complex manifold $Y$ with identically vanishing
    Kobayashi pseudo distance $d_Y$ whose infinitesimal
    Kobayashi-Royden pseudo metric is not vanishing ($F_Y\not\equiv 0$).
  \end{theorem}
 
 \section{Locally trivial fiber bundles}

 \begin{theorem}\label{bundle}
   There exists a locally trivial holomorphic fiber bundle
   $\rho:Y\to B$ with fiber isomorphic to the unit disc
   $\Delta=\{z\in\C:|z|<1\}$ such that
   $d_Y\equiv 0$.
 \end{theorem}

 \begin{proof}
   The starting point is the construction of Theorem~\ref{main}.
   Recall that $j_2$ denoted an embedding of the free group $F_2$ into
   $SL_2(\C)$ as a discrete subgroup.
   Hence the quotient $B=j_2(F_2)\backslash SL_2(\C)$ is a complex manifold.
   The natural  projection of $X=\H\times SL_2(\C)$ onto its second factor,
   $SL_2(\C)$, induces a projection from $X/\mu=Y$ onto $B$.
   By construction this is a locally trivial holomorphic fiber bundle
   with $\H\simeq\Delta$ as fiber.
   This yields the assertion of the proposition, since $d_Y\equiv 0$
   by Theorem~\ref{main}.
 \end{proof}
  
\section{Finite coverings}  
  
\begin{proposition}\label{finite}
  Let $\pi:X\to Y$ be a finite unramified covering

  Then the Kobayashi pseudodistance $d_X$ vanishes if and only if
  $d_Y$ vanishes.
\end{proposition}

\begin{proof}
  Assume $d_Y\equiv 0$.
  
  Since the covering is finite, the infimum in
  \eqref{cover-formel} is actually a minimum.
  It follows that for every $\tilde p\in X$ and $q\in Y$ there exists
  an element $\tilde q\in\pi^{-1}(q)$ with $d_X(\tilde p,\tilde q)=0$.
    As a consequence, for every $\tilde q\in Y$ we have
    $\pi(E_{\tilde p})=Y$ if we define
    \[
    E_{\tilde p}\stackrel{def}{=}\{x\in X:d_X(x,\tilde p)=0\}.
    \]
    These sets $E_{\tilde q}$ are the equivalence classes for the
    natural equivalence relation defined by $x\sim y\iff d_X(x,y)=0$.
    Hence for $x,y\in X$ either $E_x=E_y$ or $E_x\cap E_Y$ is empty.
    It follows that the set $E_x$ define a partition of $X$ into finitely
    many disjoint closed subsets. But $X$ is connected, therefore
    there exists no non-trivial finite partition
    of $X$ into disjoint closed subsets.
    Thus $E_x=X$ for any $x$, which implies $d_X\equiv 0$.
\end{proof}

\bibliographystyle{line}
\bibliography{cover}
 
\end{document}